\documentclass{amsart}
\pdfoutput=1
\usepackage{caption}
\usepackage{subcaption}
\usepackage{amssymb}
\usepackage{graphicx}
\usepackage{wrapfig}
\begin{document}

\author{G\'abor Fejes T\'oth}
\address{Alfr\'ed R\'enyi Institute of Mathematics,
Re\'altanoda u. 13-15., H-1053, Budapest, Hungary}
\email{gfejes@renyi.hu}

\title{Packing and covering properties of sequences of convex bodies}
\thanks{The English translation of the book ``Lagerungen in der Ebene,
auf der Kugel und im Raum" by L\'aszl\'o Fejes T\'oth will be
published by Springer in the book series Grundlehren der
mathematischen Wissenschaften under the title
``Lagerungen---Arrangements in the Plane, on the Sphere and
in Space". Besides detailed notes to the original text the
English edition contains eight self-contained new chapters
surveying topics related to the subject of the book but not
contained in it. This is a preprint of one of the new chapters.}

\begin{abstract}
This paper suveys different variants of the following problem:
Given a convex set $K$ and a sequence $\{C_i\}$ of convex bodies in $E^n$,
is it possible to pack the sequence of bodies in $K$ or cover $K$ with
the bodies? Algorithmic versions of these problems are on-line packing and
on-line covering: The bodies of the sequence are given one at a time and
the algorithm is to decide on the placement of the arriving body before
the next body is revealed; once placed, the body cannot be moved.
\end{abstract}

\maketitle

Consider the following problem: Given a convex set $K$ and a sequence
$\{C_i\}$ of convex bodies in $E^n$, is it possible to pack the sequence
of bodies in $K$ or cover $K$ with the bodies? More specifically, is it
possible to find for each $i$ a congruent copy $\overline{C}_i$ of each
$C_i$, so that the bodies $\{\overline{C}_i\}$ form a packing in, or a
covering of, $K$? If so, then we say that the sequence $\{C_i\}$ permits
an {\it isometric packing in} $K$, or an {\it isometric covering of} $K$,
respectively. If translates of the bodies $C_i$ are used, then we say
that $\{C_i\}$ permits a {\it translative packing in} $K$, or a
{\it translative covering of} $K$.

Problem 10.1 from the Scottish Book (see {\sc{Mauldin}} \cite[p.~74]{Mauldin})
reads as follows:

\smallskip
``PROBLEM 10.1: MAZUR, AUERBACH, ULAM, BANACH

Theorem. If $\{K_n\}_{n=1}^\infty$ is a sequence of convex bodies, each of
diameter $\le a$ and the sum of their volumes is $\le b$, then there exists
a cube with the diameter $c=f(a,b)$ such that one can put all the given bodies
in it disjointly.

Corollary. One kilogram of potatoes can be put into a finite sack.
Determine the function $c=f(a,b)$."
\smallskip

Because of the corollary, packing problems of this type are sometimes called
{\it potato-sack problems}.

The Scottish Book does not give a proof of the theorem. A proof was later
described by {\sc{Kosi\'{n}ski}} \cite{Kosinski}, who presented
an explicit bound on $f(a, b)$. The main idea is to
enclose each $C_i$ in a {\it box}, that is in a rectangular parallelepiped,
whose volume is greater than that of $C_i$ at most by a constant factor
independent from $C_i$, and whose diameter is not too much larger than the
diameter of $C_i$, thereby reducing the problem to packing a sequence of boxes
in a box. Then it is shown that every sequence of $n$-dimensional boxes of
edges at most $D$ and total volume at most $V$ can be packed in a box whose
$n-1$ edges are of length $3D$ and the $n$-th one is of length
$(V+D^n)/D^{n-1}$. {\sc{Moon}} and {\sc{Moser}} \cite{MoonMoser} improved the
above bound. They proved that such a family of boxes can be packed in a box of
sides $2D, 2D,\ldots, 2D, 2(V+D^n)/D^{n-1}$.

Moon and Moser also addressed the problem of covering and proved that a
family of $n$-dimensional boxes of edges at most $D$ and total volume $V$
can cover a cube of side $D$ if $V\ge c_n(2D)^n$, where
$c_n=\frac{2\cdot4\cdot8\cdots2^n}{1\cdot3\cdot7\cdots(2^n-1)}<2.463$.
{\sc{Meir}} and {\sc{Moser}} \cite{MeirMoser} slightly improved this result.
In particular, they showed that any family of $n$-dimensional cubes of total
volume $V$ can transitively cover a cube of volume $V/(2^n-1)$. Concerning
packing cubes in a cube they proved that a family of cubes of total volume $V$
can be transitively packed in a cube of volume $2^{n-1}V$ provided that the
side-lengths of the cubes do not exceed $2^{(n-1)/n}\sqrt[n]{V}$. The bound
$2^{n-1}V$ is best possible, since two cubes of volume $V/2$ cannot be
transitively packed in a cube of volume smaller than $2^{n-1}V$.

Let $\mathcal{F}$ be a family of convex bodies. For a given convex body $C$ we
define $p_i(C|{\mathcal{F}})$ as the greatest number such that every sequence
of members of ${\mathcal{F}}$ whose diameter does not exceed the diameter
of $C$ and whose total volume is at most $p_i(C|{\mathcal{F}})\,{\mathrm{vol}}(C)$,
permits an isometric packing in $C$. Further, we define $c_i(C|{\mathcal{F}})$
as the smallest number such that every sequence of members of $\mathcal{F}$
with total volume at least $c_i(C|{\mathcal{F}})\,{\mathrm{vol}}(C)$ permits
an isometric covering of $C$. We define similarly $p_t(C|{\mathcal{F}})$ and
$c_t(C|{\mathcal{F}})$ where the extreme values of the total volume are taken
for translative packings and translative coverings, respectively.

In most of the investigated cases
$\mathcal{F}$ consists of homothetic or similar copies of a convex body. For
a given convex body $C$ let ${\mathcal{C}}_h$ and ${\mathcal{C}}_s$ denote the family
of (positive and negative) homothetic copies of $C$ and the family of similar copies
of $C$, respectively. Further, let ${\mathcal{C}}_h^+$ and ${\mathcal{C}}_h^-$ be
the families of positive homothetic copies and negative homothetic copies
of $C$, respectively.

\subsection{Packing and covering cubes and boxes}

Let $I^n(s)$ denote an $n$-dimen\-sional cube of edge length $s$. If the size is
irrelevant, we simply write $I^n$. Then the theorems of {\sc{Meir}} and {\sc{Moser}}
\cite{MeirMoser} mentioned above state that
$$p_t(I^n|{\mathcal{I}}^n_h)=2(1/2)^n$$
and
$$c_t(I^n|{\mathcal{I}}^n_h)=2^n-1.$$

The latter equality was independently proved also by {\sc{A. Bezdek}}
and {\sc{K. Bez\-dek}} \cite{BezdekABezdekK84}.

Better bounds have been obtained under the assumption that the sets used
for packing and covering are uniformly bounded and the container is large.
Let ${\mathcal{B}}^n$ denote the family of $n$-dimensional boxes of edge
length at most 1. {\sc{Groemer}} \cite{Groemer82} proved that for $s\ge3$
$$p_i(I^n(s)|{\mathcal{B}}^n)\ge(s-1)^n-\frac{s-1}{s-2}((s-1)^{n-2}-1)$$
and
$$c_i(I^n(s)|{\mathcal{B}}^n)\le(s+1)^n-1.$$

\section{Results for general convex bodies}

{\sc{Macbeath}} \cite{Macbeath} proved that to every $n$-dimensional
convex body $C$ there exist two boxes, $Q_1$ and $Q_2$, with
$$n^nV(Q_1)\ge V(C)\ge\frac{1}{n!}V(C),$$
such that $Q_1\subset C\subset Q_2$. Let ${\mathcal{F}}^n$ denote the family
of $n$-dimensional convex bodies whose diameters are at most 1. The
combination of Macbeath's theorem with the above mentioned theorems of
Groemer yields
$$p_i(I^n(s)|{\mathcal{F}}^n)\ge
\frac{1}{n!}((s-1)^n-\frac{s-1}{s-2}((s-1)^{n-2}-1))$$
and
$$c_i(I^n(s)|{\mathcal{F}}^n)\le n^n((s+1)^n-1).$$

Macebeath's theorem was also proved by {\sc{Hadwiger}} \cite{Hadwiger55b}.
The part about the box containing $C$ was proved also by
{\sc{Kosi\'{n}ski}} \cite{Kosinski} and the part about the box contained in
$C$ by {\sc{Chakerian}} \cite{Chakerian75}. {\sc{Lassak}} \cite{Lassak93}
improved the bound for the volume of the box $Q_1$ by a factor of $2$. In
particular, every convex disk $C$ of area $a$ in the plane is contained in
a rectangle of area $2a$. This was previously proved also by {\sc{Radziszewski}}
\cite{Radziszewski} and {\sc{S\"uss}} \cite{Suss}.

Since two copies of $(\frac{1}{2}+\varepsilon)C$ cannot be packed transitively
in $C$ we have $p_t(C|{\mathcal{C}}_h^+)\le\frac{1}{2}$. {\sc{Soifer}}
\cite{Soifer99a} and {\sc{Novotny}} \cite{Novotny} conjectured that
$p_t(C|{\mathcal{C}}_h^+)=\frac{1}{2}$ for every convex disk $C$.
We are far from being able to prove or disprove this conjecture. The best
known lower bound for general convex disks is $p_t(C|{\mathcal{C}}_h^+)\ge\frac{1}{4}$
due to {\sc{Januszewski}} \cite{Januszewski07a}. Concerning packing
positive and negative homothetic copies in a convex disk {\sc{Januszewski}}
\cite{Januszewski08a}  established the inequality
$p_t(C|{\mathcal{C}}_h)\ge\frac{7}{40}$. The stronger conjecture of
{\sc{Soifer}} \cite{Soifer99a} that we have even
$p_i(C|{\mathcal{C}}_h)=\frac{1}{2}$ for every convex disk $C$ was refuted
by {\sc{Novotny}} \cite{Novotny} by showing that for the rectangle $R$ with
sides $3^{1/4}$ and $2^{1/4}$ $p_i(R|{\mathcal{R}}_h)=\sqrt{3/8}$.

L. Fejes T\'oth conjectured (see {\sc{Brass, Moser}} and {\sc{Pach}}
\cite[p.~131]{BrassMoserPach}) that $C_t(C|{\mathcal{C}}_h^+)\le3$ for every
convex disk $C$. The first upper bound of $C_t(C|{\mathcal{C}}_h^+)\le12$
for general convex disks $C$ is due to {\sc{A.~Bezdek}} and {\sc{K.~Bezdek}}
\cite{BezdekABezdekK84}. {\sc{B\'alint B\'alintov\'a, Branick\'a,
Gre\v{s}\'ak, Hrinko, Novotn\'y and Stacho}} \cite{Balint+} lowered this
bound to $9(9-\sqrt{17})/4$. Further improvements were achieved by
{\sc{Janu\-szewski}} \cite{Januszewski98a,Januszewski01,Januszewski03a}
bringing down the upper bound to $6.5$. In \cite{Januszewski98a}
{\sc{Janu\-szewski}} considered the $n$-dimensional case as well, and showed
that $c_t(C|{\mathcal{C}}_h^+)\le(n+1)^n-1$ for every convex body $C$ in $E^n$.
{\sc{Nasz\'odi}} \cite{Naszodi10} improved this bound to $6^n$ for general
convex bodies and to $3^n$ if $C$ is centrally symmetric. A further improvement
was achieved recently by {\sc{Livshyts}} and {\sc{Tikhomirov}} \cite{LivshytsTikhomirov20b}:
$c_t(C|{\mathcal{C}}_h^+)\le2^nn\log{n}(1+o(n))$ if $C$ is centrally symmetric, and
$c_t(C|{\mathcal{C}}_h^+)\le\frac{1}{\sqrt{\pi{n}}}4^nn\log{n}(1+o(n))$ otherwise.

{\sc{Soltan}} (see {\sc{Brass}}, {\sc{Moser}} and {\sc{Pach}}
\cite[Sect. 3.2, Conjecture 2]{BrassMoserPach}) formulated the following conjecture :
If $C$ is a convex body in $E^n$ and $C\subset\cup\lambda_i C_i$ for some positive
coefficients $\lambda_i<1$ then $\sum\lambda_i\ge n$. {\sc{Soltan}} and {\sc{V\'as\'arhelyi}}
\cite{SoltanVasarhelyi} proved the conjecture for $n=2$ and showed that equality
characterizes parallelograms. The special case of the conjecture when $C$ is a
triangle or a parallelogram was also proved by {\sc{Dumitrescu}} and {\sc{Jiang}}
\cite{DumitrescuJiang08}. Soltan and V\'as\'arhelyi settled the conjecture also for
the case when the number of copies covering $C$ is at most $n+1$. {\sc{Nasz\'odi}}
\cite{Naszodi10} proved the following asymptotic version of the conjecture:
For any $\nu<1$ there is an $n_0$ such that if $n>n_0$ then for every
$n$-dimensional convex body $C$, if some homothetic copies of $C$ of ratios
$0<\lambda_1, \lambda_1,\ldots,\lambda_m<1$ cover $C$ then
$\sum_{i=1}^m\lambda_\ge\nu n$. {\sc{Glazyrin}} \cite{Glazyrin19} proved
the conjecture for balls.

{\sc{Ambrus}} \cite{Ambrus2022} proved that if an
$n$-dimensional simplex is covered by its negative homothetic copies then
the sum of the absolute values of the coefficients is at least $n$.
Suggested by this result we raise the question: Is it true that if
$C$ is a convex body in $E^n$ and $C\subset\cup\lambda_i C_i$ for some
coefficients $-1<\lambda_i<0$ then $\sum|\lambda_i|\ge n$. Since, in general,
it is more difficult to cover a convex body by its negative than by its
positive homothetic copies, we risk the conjecture that the answer is
affirmative. Maybe, the same conclusion holds if we only suppose
that the absolute value of the coefficients is less than 1.

\section{On-line packing and covering}

The problem of packing a container with a sequence of convex
bodies has an algorithmic version: The bodies of the sequence are given one at
a time, as on a conveyor belt, and the algorithm is to decide on the placement
of the arriving body before the next body is revealed; once placed, the body
cannot be moved. We call this an {\it on-line} packing problem.
On-line covering problems are defined similarly. Research in this direction was
initiated by {\sc{Lassak}} and {\sc{Zhang}} \cite{LassakZhang} for packing and
by {\sc{W.~Kuperberg}} \cite{Kuperberg94a} for covering. Analogously to the
quantities $p_i(C|{\mathcal{F}})$, $c_i(C|{\mathcal{F}})$, $p_t(C|{\mathcal{F}})$
and $c_t(C|{\mathcal{F}})$ we define $p^*_i(C|{\mathcal{F}})$,
$c^*_i(C|{\mathcal{F}})$, $p^*_t(C|{\mathcal{F}})$ and $c^*_t(C|{\mathcal{F}})$
where the extreme values are taken for on-line arrangements only.

Recall that $p_t(I^n|{\mathcal{I}}^n_h)=2(\frac{1}{2})^n$. Improving on an
earlier result by {\sc{Lassak}} \cite{Lassak97c}, {\sc{Januszewski}} and
{\sc{Lassak}} \cite{JanuszewskiLassak97} proved that for $n\ge5$ an equally efficient
on-line algorithm exists:
$$p^*_t(I^n|{\mathcal{I}}^n_h)=2({1}/{2})^n.$$
For $n=3$ and $n=4$ they proved the somewhat weaker result
$p^*_t(I^n|{\mathcal{I}}^n_h)\ge\frac{3}{2}(\frac{1}{2})^n$. The 4-dimensional
case was settled recently by {\sc{Januszewski}} and {\sc{Zielonka}}
\cite{JanuszewskiZielonka18a}: $p^*_t(I^4|{\mathcal{I}}^4_h)=1/8$.
The question whether $p^*_t(I^n|{\mathcal{I}}^n_h)=2(\frac{1}{2})^n$
for $n=2$ and $n=3$ remains open.

Concerning packing boxes in a cube, {\sc{Lassak}} \cite{Lassak97c} proved that
$p^*_t(I^n|{{\mathcal{B}}_n}_h)\ge(1-\sqrt3/2)^{n-1}$, which was improved by
{\sc{Januszewski}} and {\sc{Zielonka}} \cite{JanuszewskiZielonka20} to
$p^*_t(I^n|{{\mathcal{B}}_n}_h)\ge(3-2\sqrt2)3^{-n}$. It is an open question
whether $p^*_t(I^n|{{\mathcal{B}}_n}_h)=p_t(I^n|{{\mathcal{B}}_n}_h)=2(\frac{1}{2})^n$
for $n\ge2$.

The algorithm by {\sc{W.~Kuperberg}} \cite{Kuperberg94a} yields
$c^*_t(I^n|{\mathcal{I}}^n_h)\le4^n$. Better algorithms by {\sc{Januszewski}}
and {\sc{Lassak}} \cite{JanuszewskiLassak94} and by {\sc{Lassak}} \cite{Lassak95}
provide the bound $c^*_t(I^n|{\mathcal{I}}^n_h)\le2^n(1+o(n))$. The breakthrough
was achieved by {\sc{Januszewski, Lassak, Rote}} and {\sc{Woeginger}}
\cite{JanuszewskiLassakRoteWoeginger96a}, who constructed an on-line
algorithm showing that $c^*_t(I^n|{\mathcal{I}}^n_h)\le2^n+3-\frac{2^{n+2}-8}{2^{2n}+2^-2}$.
This bound comes astoundingly close to the best value for off-line coverings.
{\sc{Lassak}} \cite{Lassak02} further improved this bound to
$$c^*_t(I^n|{\mathcal{I}}^n_h)\le2^n+\frac{5}{3}(1 + 2^{-n}).$$
For the 3-dimensional case this yields $c^*_t(I^3|{\mathcal{I}}^3_h)\le8+15/8=9.875$.

The main tool used in the two articles cited above is the
{\it{$q$-adic on-line algorithm}} for covering the interval $[0,1]$
with a sequence of segments $S_i$ of length $\delta_i$
$(i=1,2,\ldots)$, where $q\ge2$ is an integer,
$\delta_i\in\left\{ q^{-1}, q^{-2}, \ldots\right\}$, and each $S_i$
must be placed on one of the intervals of the form
$[k\delta_i, (k+1)\delta_i]\subset[0,1]$, for some integer $k$.
This approach was earlier suggested by W. Kuperberg \cite{Kuperberg94b},
explicitly for $q=2$ and implicitly for all $q\ge2$. The suggestion
was put in the form of a question asking for the existence of a winning
algorithm in a $2$-adic interval covering game between two players.
The solution of the problem appeared in {\sc{Januszewski, Lassak, Rote}}
and {\sc{Woeginger}} \cite{JanuszewskiLassakRoteWoeginger96b}.

\section{Special convex disks}

A considerable amount of research has been devoted to packing and covering
problems involving sequences of special convex disks, such as squares,
rectangles, or triangles.

For the special case of a square $S$ the theorem of {\sc{Meir}} and {\sc{Moser}}
\cite{MeirMoser} mentioned above states that
$p_t(S|{\mathcal{S}}_h)=\frac{1}{2}$ and $c_t(S|{\mathcal{S}}_h)=3$.
{\sc{Januszewski}} \cite{Januszewski02a} proved that
$$c_i(S|{\mathcal{S}}_h)=2$$
and in his papers \cite{Januszewski07b,Januszewski08b} he proved that
$$c_t(S|{\mathcal{S}}_s)=3.$$
Let ${\mathcal{S}}'$ be the family of squares with diagonals parallel
to the sides of $S$. {\sc{Januszewski}} \cite{Januszewski10b,Januszewski02b}
proved that
$$p_t(S|{\mathcal{S}}')=4/9$$
and
$$c_t(S|{\mathcal{S}}')=2.5.$$
Concerning packing rectangles into rectangles {\sc{Yuan, Xu}} and {\sc{Ding}}
\cite{YuanXuRen} proved the following: If ${\mathcal{R}}_a$ denotes the
family of rectangles of side lengths not greater than $a$, and $R_{ab}$
denotes a rectangle with sides $a$ and $b\ge{a}$, then
$$p_i(R_{ab}|{\mathcal{R}}_a)=\frac{ab}{2}.$$
The case $a=b$ was proved earlier by {\sc{Januszewski}} \cite{Januszewski00}.

For on-line packing squares into a square
{\sc{Januszewski}} and {\sc{Lassak}} \cite{JanuszewskiLassak97} established
the bound $p^*_t(S|{\mathcal{S}}_h)\ge\frac{5}{16}$. The lower bound on
$p^*_t(S|{\mathcal{S}}_h)$ was subsequently improved to $\frac{1}{3}$
by {\sc{Han, Iwama}} and {\sc{Zhang}} \cite{HanIwamaZhang}, to $\frac{11}{32}$
by {\sc{Fekete}} and {\sc{Hoffmann}} \cite{FeketeHoffmann17} and, finally, to
$\frac{2}{5}$ by {\sc{Brubach}} \cite{Brubach}. Concerning on-line packing
rectangles in a square, {\sc{Lassak}} \cite{Lassak97c} proved the bound
$p^*_i(S|{\mathcal{R}})\ge\frac{5}{36}$, which was improved by {\sc{Januszewski}}
and {\sc{Zielonka}} \cite{JanuszewskiZielonka18b} to
$p^*_i(S|{\mathcal{R}})\ge0.28378$.

For on-line covering a square with squares {\sc{Januszewski}} and {\sc{Lassak}}
\cite{JanuszewskiLassak95a} proved
$c^*_t(S|{\mathcal{S}}_h)\le\frac{7}{4}\sqrt{3}{9} + \frac{13}{8}\approx 5.265$.
This was improved to $c^*_t(S|{\mathcal{S}}_h)\le4$ by {\sc{Januszewski}}
\cite{Januszewski09a}. Concerning on-line covering a square by rectangles
{\sc{Januszewski}} \cite{Januszewski96} proved
$c^*_i(S|{\mathcal{R}})\le\frac{15}{2}$.

{\sc{Richardson}} \cite{Richardson} proved that
$p_i(T|{\mathcal{T}}_s)\ge\frac{1}{2}$ for every triangle $T$. In fact,
his algorithm used only positive and negative homothetic copies of $T$
and he conjectured that the packing is possible by using positive
homothetic copies only. This conjecture is equivalent to the statement
that $p_t(T|{\mathcal{T}}_h^+)=\frac{1}{2}$. In \cite{Januszewski02c}
{\sc{Januszewski}} established the bound
$p_t(T|{\mathcal{T}}_h^+)\ge\frac{19}{56}$ and later in \cite{Januszewski09b}
he refined the method and confirmed the Richardson's conjecture. On the other hand,
{\sc{Januszewski}} disproved {\sc{Soifer}}'s conjecture \cite{Soifer99b}
that $p_i(T|{\mathcal{T}}_s)=\frac{1}{2}$ for every triangle $T$. In
\cite{Januszewski03b} he showed that $p_i(T|{\mathcal{T}}_s)=\frac{1}{2}$
if and only if $T$ is equilateral. Moreover, in \cite{Januszewski03c}
he proved that if $T$ is an isosceles right triangle then
$0.511\le\frac{1}{3}(5-1\sqrt3)\le p_i(T|{\mathcal{T}}_h)
\le\frac{7}{2}-2\sqrt2\le0.6715$. For translative packing of
positive and negative homothetic copies of a triangle {\sc{Januszewski}}
\cite{Januszewski06} proved
$$p_t(T|{\mathcal{T}}_h)=\frac{2}{9}.$$

Concerning covering a triangle $T$ with homothetic copies {\sc{V\'as\'arhelyi}}
\cite{Vasarhelyi84} proved that
$$c_t(T|{\mathcal{T}}_h^-)=4$$
and {\sc{F\"uredi}} \cite{Furedi03} proved that
$$c_t(T|{\mathcal{T}}_h^+)=2.$$
{\sc{Janu\-szewski}} \cite{Januszewski98b} strengthened V\'as\'arhelyi's result
by showing that
$$c_t(T|{\mathcal{T}}_h)=4.$$
For a right isosceles triangle $T$ {\sc{F\"uredi}} \cite{Furedi07} established
the equality
$$c_i(T|{\mathcal{T}}_h)=\frac{1+\sqrt2}{2}.$$
Denote by $T_\varphi$ the triangle obtained from $T$ by a rotation through the angle
$\varphi$. {\sc{V\'as\'arhelyi}} \cite{Vasarhelyi93} proved that
$$c_t(T_{30^\circ}|{{\mathcal{T}}_{30^\circ}}_h^+)=4$$
and $c_t(T_\varphi|{{\mathcal{T}}_\varphi}_h^+)<4$ for every triangle that is not
regular.

{\sc{Januszewski}} \cite{Januszewski09c} proved that for an equilateral triangle
$T$ and a square $S$ having a side parallel to a side of $T$ we have
$$p_t(T|{\mathcal{S}}_h)=2\sqrt3-3$$
and
$$c_t(T|{\mathcal{S}}_h)=2\sqrt3.$$
The theorem concerning covering was extended by {\sc{Lu}} and {\sc{Su}}
\cite{LuSu18} to covering an isosceles triangle $T(h)$ with base length
$1$ and with height $h$ by homothetic copies of a square $S$ one side
of which is parallel to the base of $T(h)$. They showed that
\[
c_t(T(h)|{\mathcal{S}}_h) =
  \begin{cases}
    \frac{2}{h} & \text{if } \frac{\sqrt2}{2}\le{h}\le1,\\
    4h & \text{if } \frac{\sqrt2}{2}\le{h}\le1,\\
    \frac{4}{h} & \text{if } 1\le{h}<\sqrt2,\\
    2 & \text{if } \sqrt2\le{h}.
  \end{cases}
\]
For a right triangle $T_0$ with
legs 1 and $\sqrt2$ and a square $S$ with sides parallel to the legs of $T$
{\sc{Fu, Lian}} and {\sc{Zhang}} \cite{FuLianZang} proved the inequality
$p_t(T_0|{\mathcal{S}}_h)\ge\frac{16-6\sqrt2}{23}$. These authors also investigate
the problem of covering a tetrahedron $T_r$ with three mutually perpendicular
edges of lengths 1, 1, and $\sqrt2$ by homothetic copies of a cube $C$ with sides
parallel to the edges of $T_r$. They prove that $c_t(T_r|{\mathcal{C}}_h)\le6\sqrt2+1$.

Concerning packing circles in a circle {\sc{Fekete, Keldenich}} and {\sc{Scheffer}}
\cite{FeketeKeldenichScheffer} proved that
$$p_t(B^2|{\mathcal{B}}^2_h)=\frac{1}{2}.$$

For the on-line case  {\sc{Januszewski}} \cite{Januszewski11} established the bound
$p^*_t(B^2|{\mathcal{B}}^2_h)>0.197$. For the corresponding covering problem
{\sc{Dumitrescu}} and {\sc{Jiang}} \cite{DumitrescuJiang10} proved that
$c_t(B^2|{\mathcal{B}}^2_h)\le2.97$ thereby confirming the conjecture of L. Fejes
T\'oth for the first convex disk that is not a polygon. They also considered the
on-line version of the problem and proved that $c^*_t(B^2|{\mathcal{B}}^2_h)\le9.763$.
{\sc{Januszewski}} \cite{Januszewski11} improved the latter bound to
$c^*_t(B^2|{\mathcal{B}}^2_h)<6.488$.

{\sc{Fekete, Morr}} and {\sc{Scheffer}} \cite{FeketeMorrScheffer} investigated
the problem of packing sequences of circles in a square or triangle. They proved
that $$p_t(S|{\mathcal{B}}^2_h)=\frac{\pi}{3+2\sqrt2}.$$
Further, if $T$ is a non-acute triangle with an incircle of area $a$ then
$$p_t(T|{\mathcal{B}}^2_h)=a.$$
For on-line packing circles in a square {\sc{Fekete, von H\"oveling}} and {\sc{Scheffer}}
\cite{FeketeHovelingScheffer} proved the inequality $p_t^*(S|{\mathcal{B}}^2_h)\le0.350389$.
For packing squares in a circle {\sc{Fekete, Gurunathan, Juneja, Keldenich, Kleist,}}
and {\sc{Scheffer}} \cite{FeketeGurunathanJunejaKeldenichKleistScheffer} established
the equality
$$p_t(B^2|{\mathcal{S}}_h)=\frac{8}{5\pi}.$$

The problem of covering the square by a sequence of circular disks was solved
by {\sc{Fekete, Gupta, Keldenich, Scheffer}} and {\sc{Shah}}
\cite{FeketeGuptaKeldenichSchefferShah}. They proved that
$$c_t(S|{\mathcal{B}}^2_h)=\frac{195\pi}{256}.$$
More generally, they gave an algorithm for covering the rectangle $R(1,\lambda)$
with sides $1$ and $\lambda\ge1$, and showed that there is a threshold value
$\lambda_0=\sqrt{\sqrt{7}/2-1/4}=1.035797\ldots$, such that for $\lambda<\lambda_0$
$$c_t(R(1,\lambda)|{\mathcal{B}}^2_h)=3\pi\left(\frac{\lambda^2}{16}+\frac{5}{32}+\frac{9}{256\lambda^2}\right),$$
and for $\lambda\ge\lambda_0$
$$c_t(R(1,\lambda)|{\mathcal{B}}^2_h)=\frac{(\lambda^2+2)\pi}{4}.$$

\section{Packing in and covering of the whole space}

The investigation of covering the whole space by sequences of convex
bodies was initiated by {\sc{Chakerian}} \cite{Chakerian75}. Clearly,
for a sequence $\{C_i\}$ of convex bodies to permit a covering of
space it is necessary that the sum of the volumes $V(C_i)$, as well as the sum
of the widths $w(C_i)$ be divergent. These conditions are not sufficient.
{\sc{Chakerian}} and {\sc{Groemer}} \cite{ChakerianGroemer74} gave necessary
and sufficient conditions for a sequence of convex disks to permit a covering
of the plane. A sequence $\{C_i\}$ of convex disks permits a covering of the
plane if and only if either the total area of the subsequence with diameter
at most 1 is infinite or the total width of the subsequence with diameter
greater than 1 is infinite. The sequence $\{C_i\}$ of convex bodies is
{\it bounded} if the set of the diameters of the bodies is bounded. In
particular, it follows that a bounded sequence of convex disks permits a
covering of the plane if and only if the total area of the the disks
is infinite. The analogous statement for $E^n$, $n\ge3$ was proved by
{\sc{Groemer}} \cite{Groemer76}. {\sc{Chakerian}} and {\sc{Groemer}}
\cite{ChakerianGroemer78} gave necessary and sufficient conditions for
a sequence of convex bodies to permit a covering of $E^n$. {\sc{Groemer}}
\cite{Groemer79} proved that for a sequence of convex bodies to
permit a covering of $E^n$ it is necessary and sufficient that the
sequence permits a translative covering of almost all points of $E^n$.
{\sc{Groemer}} \cite{Groemer80,Groemer83a} investigated coverings of
space by finite sequences of unbounded convex sets and in
\cite{Groemer83b} he gave an upper bound for the total volume of
a sequence of convex bodies permitting a covering of the $n$-dimensional
sphere.

As consequences of the results of {\sc{Groemer}} \cite{Groemer82} mentioned
in Section 1.1, we note the following. If $\{C_i\}$ is a bounded sequence of
$n$-dimensional convex bodies such that $\sum V(C_i)=\infty$, then it permits
an isometric covering of the $n$-dimensional space with density $\frac{1}{2}n^n$
and an isometric packing with density $\frac{1}{n!}$. Moreover, if all the
sets $C_i$ are boxes, then $\{C_i\}$ permits a translative covering of space
and a translative packing in space with density~$1$.
In the plane, any bounded sequence $\{C_i\}$ of convex disks with $\sum
C_i=\infty$ permits even a translative packing and covering with density
$\frac{1}{2}$ and $2$, respectively. It is an open problem whether for $n>2$
every bounded sequence $\{C_i\}$ of $n$-dimensional convex bodies of infinite
total volume permits a translative covering of $E^n$.

{\sc{Groemer}} \cite{Groemer88} investigated the question under which conditions
a sequence of convex bodies $\{C_i\}$ in $E^n$ permits an isometric packing
or covering of density 1. He showed that the conditions $\sum V(C_i)$ and
$\lim_{i\to\infty}d(C_i)=0$ on the volume and diameter of the sets $C_i$
are sufficient for constructing such packings and coverings.

\section{Covering with slabs}

Concerning the problem of covering space with a sequence of slabs,
{\sc{Groemer}} \cite{Groemer81a, Groemer81b} proved that every sequence of
slabs of widths $w_i$ in $E^n$ for which $\sum w_i^{(n+1)/2}=\infty$ permits a
translative covering. {\sc{Makai}} and {\sc{Pach}} \cite{MakaiPach}
conjectured that a sequence of slabs in $E^n$ permits a translative covering
if and only if the sum of their widths is infinite. They verified the
conjecture for the case of the plane. In higher dimensions the conjecture is
still unresolved. {\sc{Kupavskii}} and {\sc{Pach}} \cite{KupavskiiPach} proved
that if $w_1\ge w_2\ge\ldots$ is an infinite sequence of positive numbers
such that
$$
\limsup_{i\to\infty}\frac{w_1+w_2+\ldots w_i}{\ln(1/w_i)}>0,
$$
then every sequence of slabs of width $w_i$ ($i=1,2,\ldots$) permits a
translative covering of $E^n$. With this result they got rather close to the
proof of the conjecture: For example, it implies that a sequence of slabs
of width $w_i=1/i$ ($i=1,2,\ldots$) permits a translative covering, while
this is false for the sequence of widths $w_i=1/i^{1+\varepsilon}$ for
any $\varepsilon>0$.

A detailed account on the topic of packing and covering properties of
sequences of convex bodies is found in the survey by {\sc{Groemer}}
\cite{Groemer85}. For surveys on on-line algorithms see {\sc{Lassak}}
\cite{Lassak97a} and {\sc{Csirik}} and {\sc{Woeginger}} \cite{CsirikWoeginger}.

\small{
\bibliography{pack}}